\documentclass[12pt]{article}

\usepackage{lineno,hyperref}
\usepackage{amsmath,amssymb,soul}
\usepackage{graphicx}

\newcommand{\ind}{\operatorname{ind}}
\newcommand{\norm}[1]{\left\Vert#1\right\Vert}
\newcommand{\scal}[1]{\langle#1\rangle}
\newcommand{\abs}[1]{\left\vert#1\right\vert}
\newcommand{\eq} [1] {\begin{equation}\label{#1}\quad}
\newcommand{\en} {\end{equation}}

\newcommand{\C}{{\mathbb C}} \newcommand{\R}{{\mathbb R}}
\newcommand{\diag}{\operatorname{diag}}
\newcommand{\dist}{\operatorname{dist}}
\newcommand{\tr}{\operatorname{trace}}
\newcommand{\im}{\operatorname{Im}}
\newcommand{\Ker}{\operatorname{Ker}}

\begin{document}

\begin{center}
{\LARGE \bf Robert Sheckley's Answerer\\[0.2ex]for two orthogonal projections}

\vspace{10mm}
{\Large Albrecht B\"ottcher and Ilya Spitkovsky\footnote{This author was supported 
in part by Faculty Research funding from the Division of Science and Mathematics, New York University Abu Dhabi}}
\end{center}

\vspace{7mm}
\begin{quote}
{\bf Abstract.} 
This paper is the written version of our talk
(presented by the second author) at the IWOTA in Chemnitz in August 2017.
The meta theorem of the paper is that Halmos' two projections theorem
is something like Robert Sheckley's Answerer:
no question about the $W^*$- and $C^*$-algebras generated by two orthogonal projections will go
unanswered, provided the question is not foolish. An alternative approach to
questions about two orthogonal projections makes use of the supersymmetry equality introduced
by Avron, Seiler, and Simon. A noteworthy
insight of the paper reveals that the supersymmetric approach is nothing but Halmos in different language
and hence an equivalent Answerer.


\end{quote}

\section{Introduction} \label{s:intro}
One of the books which had a great influence on us when we just started studying
Functional Analysis was Glazman and Lyubich's~\cite{GlaLyu69}. In particular, we always remembered Glazman's famous
``And how does this look in the two-dimensional case?'' question when someone was describing to
him an elaborate infinite-dimensional construction, and the claim that ``quite frequently this shocking question helped to
better understand the gist of the matter''. The topic of this paper is a striking example of the validity of Glazman's approach.

So, let us start with a pair of orthogonal projections $P,Q$ acting on $\C^2$.
If one of them, say $P$, is the zero or the identity operator, we may diagonalize $Q$ by a unitary similarity to
$\diag[0,0]$, $\diag[1,0]$, or $\diag[1,1]$, while $P$ remains equal to $\diag[0,0]$ or $\diag[1,1]$
under this unitary similarity. Thus suppose $P,Q$ both have rank one.
A unitary similarity
can then be used to put $P$ in the diagonal form $\diag[1,0]$. The matrix of $Q$ in the
respective basis is Hermitian, with zero determinant and the trace equal to one. An additional (diagonal) unitary similarity,
while leaving the representation of $P$ unchanged, allows us to make the off-diagonal entries of this matrix equal and non-negative,
without changing its diagonal entries. It is thus bound to equal
\[ \begin{bmatrix} x & \sqrt{x(1-x)} \\  \sqrt{x(1-x)} & 1-x \end{bmatrix} \]
with $x\in (0,1)$ (the values $x=0,1$ are excluded because otherwise $Q$ would commute with $P$).

\section{Canonical representation} \label{s:canrep}

This picture extends to the general Hilbert space setting in the most natural and direct way.
Namely, according to Halmos' paper~\cite{Hal69}, for a
pair of orthogonal projections acting on a Hilbert space $\mathcal H$ there exists an orthogonal decomposition
\eq{deco} {\mathcal H}={\mathcal M}_{00}\oplus{\mathcal M}_{01}\oplus{\mathcal M}_{10}
\oplus{\mathcal M}_{11}\oplus\left({\mathcal M}\oplus{\mathcal M}'\right), \en
with respect to which
\eq{canrep} \begin{aligned} P &= & I\oplus I\oplus 0\oplus 0 \;\oplus & \: \boldsymbol{W}^*
\begin{bmatrix}I & 0 \\ 0 & 0\end{bmatrix}\boldsymbol{W} , \\
Q & = &  I\oplus 0\oplus I\oplus 0 \;\oplus & \;\boldsymbol{W}^*\begin{bmatrix} H & \sqrt{H(I-H)}\; \\
\sqrt{H(I-H)} & I-H \; \end{bmatrix}\boldsymbol{W},\end{aligned} \en
where $\boldsymbol{W}=\begin{bmatrix} I & 0 \\ 0 & W\end{bmatrix}$,
$W\colon{\mathcal M}'\longrightarrow{\mathcal M}$ is unitary, and $H$ is
the compression of $Q$  to ${\mathcal M}$. The operator $H$
is selfadjoint with spectrum $\sigma(H)\subset[0,1]$ and $0,1$ not being its eigenvalues.
We refer to~\cite{BSpit10, Roch} for more on the history of this representation before and after Halmos,
for full proofs, and for related topics.
One more proof will be given in Section \ref{s: super}.

Of course,
\begin{eqnarray*}
& & {\mathcal M}_{00}=\im P\cap\im Q, \quad {\mathcal M}_{01}=\im P\cap\Ker Q,\\
& & {\mathcal M}_{10}=\Ker P\cap\im Q,\quad {\mathcal M}_{11}=\Ker P\cap\Ker Q,
\end{eqnarray*}
and so
\eq{M} {\mathcal M}=\im P\ominus\left({\mathcal M}_{00}\oplus{\mathcal M}_{01}\right), \en
while
\[ {\mathcal M}'=\Ker P\ominus\left({\mathcal M}_{10}\oplus{\mathcal M}_{11}\right). \]
It is an implicit consequence of \eqref{canrep} that $\dim{\mathcal M}'=\dim{\mathcal M}$.
In what follows, for simplicity of notation we will
identify ${\mathcal M}'$ with ${\mathcal M}$ via their isomorphism $W$.
In other words, we will drop the factors $\boldsymbol{W}, \boldsymbol{W}^*$
in \eqref{canrep}.

The operators $P$ and $Q$ commute if and only if the last summand in~\eqref{canrep} is missing,
that is, ${\mathcal M} (={\mathcal M}')=\{0\}$. This $P,Q$ configuration
is of course not very interesting, though should be accounted for. Another extreme is ${\mathcal M}_{ij}=\{0\}$
for all $i,j=0,1$. If this is the case, $P$ and $Q$ are said
to be in the {\em generic position}.

\section{Algebras} \label{s:alg}

Based on \eqref{canrep}, a description of the von Neumann algebra ${\mathcal A}(P,Q)$ generated
by $P$ and $Q$ was obtained in \cite{GiKu}. The elements of
${\mathcal A}(P,Q)$ are all the operators of the form
\eq{elal} \left(\oplus a_{ij}I_{{\mathcal M}_{ij}}\right)\oplus\begin{bmatrix}\phi_{00}(H) & \phi_{01}(H) \\
\phi_{10}(H) & \phi_{11}(H) \end{bmatrix},  \en
where $a_{ij}\in\C$, the direct sum in the parentheses is taken with respect to $i,j=0,1$ for
which ${\mathcal M}_{ij}\neq\{0\}$ and the functions $\phi_{ij}$
are Borel-measurable and essentially (with respect to the spectral measure of $H$) bounded on $[0,1]$.

With the notation $\Phi=\begin{bmatrix}\phi_{00} & \phi_{01} \\ \phi_{10} & \phi_{11}\end{bmatrix}$,
we can (and sometimes will) abbreviate \eqref{elal} to
\eq{elal1}   \left(\oplus a_{ij}I_{{\mathcal M}_{ij}}\right)\oplus\Phi(H). \en
Invoking the spectral representation
\[H=\int_{\sigma(H)}\lambda\, dE(\lambda)\] of $H$, we can also rewrite \eqref{elal} as
\[ \left(\oplus a_{ij}I_{{\mathcal M}_{ij}}\right)\oplus\int_{\sigma(H)}\Phi(\lambda) dE(\lambda) .\]

The elements of the $C^*$-algebra ${\mathcal B}(P,Q)$ generated by $P$ and $Q$ are distinguished
among those of the form \eqref{elal} by the following \cite{Ped68,VaSpit}
additional properties\footnote{Of course, conditions on $a_{ij}$ below are meaningful
only if the respective subspaces ${\mathcal M}_{ij}$ are non-zero.}:

(i) The functions $\phi_{ij}$ are continuous on $[0,1]$, not just measurable;

(ii) If $0\in\sigma(H)$, then $\phi_{01}(0)=\phi_{10}(0)=0$, $a_{00}=\phi_{11}(0)$, $a_{11}=\phi_{00}(0)$;

(iii) If $1\in\sigma(H)$, then $\phi_{01}(1)=\phi_{10}(1)=0$, $a_{01}=\phi_{11}(1)$, $a_{10}=\phi_{00}(1)$.

In the finite dimensional setting the algebras ${\mathcal A}(P,Q)$ and ${\mathcal B}(P,Q)$ of course coincide,
and their elements are (up to a unitary similarity
which we agreed to ignore) of the form
\eq{elalf} \left(\oplus a_{ij}I_{{\mathcal M}_{ij}}\right)\oplus
\left(\bigoplus_{\lambda_j\in\sigma(H)}\Phi(\lambda_j)\right).  \en

\section{The Answerer} \label{s:ans}

Independently of whether $\mathcal H$ is finite- or infinite-dimensional,
the representations \eqref{elal}--\eqref{elalf}
allow us to settle any meaningful question about
operators from the algebras generated by the pair $P,Q$. The real challenge is to ask the right questions,
and this brings us to Robert Sheckley's famous short story ``Ask a foolish question'',
written in 1953.

In that story,
we encounter an Answerer, a machine built a long time ago by a race and left back on a planet after the race
disappeared. ``He [the Answerer] knew the nature of things, and why things are as they are, and what they
are, and what it all means.
Answerer could answer anything, provided it was a legitimate question.'' For example, he could not
give an answer to the question
``Is the universe expanding?'' What he replied was
`` `Expansion' is a term inapplicable to the situation. Universe, as the Questioner views it, is an illusory
concept.'' Another drastic passage in the story says "Imagine a bushman
walking up to a physicist and asking him why he can't shoot his arrow into the sun. The scientist can explain it
only in his own terms. What would happen?" -- "The scientist wouldn't even
attempt it, ... he would know the limitations of the
questioner.'' -- "How do you explain the earth's rotation to a bushman? Or better, how do you
explain relativity to him, maintaining scientific rigor in your explanation at all times, of course.''
-- ``We're bushmen. But the gap is much greater here. Worm and super-man, perhaps. The worm desires to know
the nature of dirt, and why there's so much of it.''
The quintessence of the story is
that

\medskip
\centerline{``In order to ask a question you must already know most of the answer.''}

\medskip
In what follows we embark on some questions about two orthogonal projections we consider as meaningful
and will show what kind of answer Halmos' theorem will give.

\section{Routine} \label{s:rou}

Some necessary bookkeeping was performed in \cite{Spit94}. An explicit, though somewhat cumbersome,
description was provided there for the kernels and ranges of operators
$A\in{\mathcal A}(P,Q)$. Based on those, Fredholmness and invertibility criteria,
formulas for spectra and essential spectra, norms, and the Moore-Penrose inverses $A^\dag$ (when it exists)
were derived.

To give a taste of these results, here is the description of $\Ker A$ for $A$ given by \eqref{elal}.
Let ${\mathcal M}^{(r)}$ be the spectral subspace of $H$ corresponding to the
subset $\Delta_r$ of $\sigma(H)$ on which $\Phi(t)$ has rank $r \in \{0,1,2\}$. Let also
\[ \phi=\sum_{i,j=0,1}\abs{\phi_{ij}}^2, \quad \chi_i=\sqrt{\frac{\abs{\phi_{0i}}^2+\abs{\phi_{1i}}^2}{\phi}}, \ i=0,1, \] and
\[ u=\exp\left(i\arg(\phi_{01}\overline{\phi_{00}}+\phi_{11}\overline{\phi_{10}})\right). \] Then
\eq{ker} \Ker A = \left(\bigoplus_{a_{ij}=0}{\mathcal{M}_{ij}}\right)\oplus\left({\mathcal M}^{(0)}\oplus{\mathcal M}^{(0)}\right)
\oplus\begin{bmatrix}u(H)\chi_1(H) \\ -\chi_0(H) \end{bmatrix}({\mathcal M}^{(1)}). \en
For example,  let $A=I-Q$. In its representation \eqref{elal} we then have
\eq{a} a_{00}=a_{10}=0, \quad  a_{01}=a_{11}=1, \en
\[\phi_{00}(t)=1-t, \quad \phi_{11}(t)=t, \quad \phi_{01}(t)=\phi_{10}(t)=-\sqrt{t(1-t)}. \]
Consequently,
\eq{fu} \phi=1, \quad u= -1, \quad \chi_1(t)=\sqrt{t}, \quad \chi_0(t)=\sqrt{1-t}. \en
Plugging \eqref{a}, \eqref{fu} into \eqref{ker}
we see that
\eq{imQ} \im Q \ (=\Ker (I-Q)) = {\mathcal M}_{00}\oplus {\mathcal M}_{10}\oplus
\begin{bmatrix}\sqrt{H}\\ \sqrt{I-H}\end{bmatrix}(\mathcal M).   \en

Since $A^*$ belongs to ${\mathcal A}(P,Q)$ along with $A$, the description of $\Ker A^*$ follows
from \eqref{ker} via a simple change of notation. The closures of $\im A^*$ and $\im A$
can then be obtained as the respective orthogonal complements. Note however that \cite{Spit94}
provides the description of these ranges themselves, not just their closures.
In particular, $\im A$ and $\im A^*$ are closed if and only if
\eq{clos} \det\Phi \text{ and } \phi \text{ are separated from 0 on }\Delta_2
\text{ and } \Delta_1 \text{ respectively,} \en
so \eqref{clos} is also a criterion for $A^\dag$ to exist. In its turn, $A$ is invertible
if and only if $\det\Phi$ is separated from zero on the whole $\Delta$
and, in addition, $a_{ij}\neq 0$ whenever ${\mathcal M}_{ij}\neq\{0\}$.

Example. Consider $A=P-Q$. Its representation \eqref{elal1} has the form
\eq{p-q} 0_{{\mathcal M}_{00}}\oplus I_{{\mathcal M}_{01}}\oplus (-I)_{{\mathcal M}_{10}}\oplus
0_{{\mathcal M}_{11}}\oplus \begin{bmatrix} I-H & -\sqrt{H(I-H)} \\
 -\sqrt{H(I-H)} & H-I\end{bmatrix}, \en and so the respective matrix $\Phi$ is
\eq{fpq} \Phi_{P-Q}(t)=\begin{bmatrix} 1-t & -\sqrt{t(1-t)} \\  -\sqrt{t(1-t)} & t-1\end{bmatrix} \en
with the characteristic polynomial $\lambda^2+t-1$. It immediately follows that
 \eq{spq} \sigma(P-Q)=\{\pm\sqrt{1-t}\colon t\in\sigma(H)\}, \en
which is a subset of $[-1,1]$ that is symmetric about the origin, with the additional
eigenvalues $1,-1$ or $0$ materializing if and only if the respective subspace ${\mathcal M}_{01}$,
${\mathcal M}_{10}$, or ${\mathcal M}_{00}\oplus{\mathcal M}_{11}$ is non-trivial.

\section{Anticommutators} \label{s:aco}

To provide yet another example of how easily the considerations of Section~\ref{s:rou}
generate some nice formulas, we turn to the anticommutator $PQ+QP$ of $P,Q$. For simplicity,
take $P$ and $Q$ in generic position. Then
\eq{aco}  PQ+QP-\lambda I=\begin{bmatrix}2H-\lambda I & \sqrt{H(I-H)} \\ \sqrt{H(I-H)}  & -\lambda I\end{bmatrix}. \en
Since the entries of the operator matrix on the right-hand side of \eqref{aco} commute pairwise,
according to \cite[Problem 70]{Hal82}
it is invertible only simultaneously with its formal determinant
\[\lambda^2 I-2\lambda H-H+H^2=(\lambda I-H)^2-H
=(\lambda I-H+\sqrt{H})(\lambda I-H-\sqrt{H}). \]
Consequently, \[ \sigma(PQ+QP)=\{\lambda\pm\sqrt{\lambda}\colon\lambda\in\sigma(H)\}. \]
In particular, $PQ+QP$ is invertible if and only if $0,1\notin\sigma(H)$.
Note that this is always the case if $\dim{\mathcal H}<\infty$.

On the other hand, in our setting $H$ is simply the operator $PQP$ considered on $\im P$,
and thus $\sigma(H)\cup\{0\}=\sigma(PQP)$. We therefore conclude that the spectrum of the anticommutator $PQ+QP$
is the set $\{\lambda\pm\sqrt{\lambda}\colon\lambda\in\sigma(PQP)\}$ from which the origin
should be removed if $0,1\notin\sigma(H)$.
This covers the result of  \cite{DouDuWang17}.

Moreover, since $PQ+QP$ is a positive semi-definite operator, its norm coincides
with the maximum of its spectrum. Therefore,
\eq{ano} \norm{PQ+QP}=\max\{\lambda+\sqrt{\lambda}\colon\lambda \in\sigma(PQP)\}
= \norm{PQP}+\norm{PQP}^{1/2}. \en
In its turn, $\norm{PQP}=\norm{PQ(PQ)^*}=\norm{PQ}^2$, and \eqref{ano} can be rewritten as
\[ \norm{PQ+QP}=\norm{PQ}^2+\norm{PQ}. \]
The latter formula was the main subject of Walters' \cite{Walt16}.

\section{Drazin invertibility} \label{s:dra}

Recall that an operator $A$ acting on a Hilbert (or even a Banach) space is {\em Drazin invertible}
if and only if the sequences $\im A^j$ and $\Ker A^j$ stabilize. If this is the case, and $k$
is the smallest non-negative integer for which $\Ker A^k=\Ker A^{k+1}$ and $\im A^k=\im A^{k+1}$,
the {\em Drazin inverse} $X$ of $A$ is defined uniquely by the properties
\[ A^{k+1}X=A^k, \quad XAX=X, \quad AX=XA. \]
A criterion for Drazin invertibility of operators $A\in{\mathcal A}(P,Q)$ and a
formula for their Drazin inverse $A^D$ was found in \cite{BSpit09}.
Setting $\Delta_{11}:=\{t\in\Delta_1\colon\tr\Phi(t)\neq 0\}$,
we have that $A$ is Drazin invertible if and only if
\eq{dra}  \det\Phi| \Delta_2 \;\:\mbox{and}\;\: \tr\Phi|\Delta_{11} \;\:\mbox{are separated from 0.}\en
Note that the first parts of conditions \eqref{clos}, \eqref{dra} are the same, while the
second requirement of \eqref{dra} implies that $\phi$ is separated from zero on $\Delta_{11}$
though not necessarily on the whole $\Delta_1$. So, if $\Delta_{10}:=\Delta_1\setminus\Delta_{11}\neq\emptyset$,
a Drazin invertible operator $A$ may or may not have closed range and thus
be Moore-Penrose invertible or not (and, even if it is, $A^D\neq A^\dag$).
This is exactly the case when $k=2$. On the other hand, if $\Delta_{10}=\emptyset$,
then condition \eqref{dra} implies \eqref{clos}. So, $A$ is Moore-Penrose invertible with
$A^D=A^\dag$ and $k$ is either zero (in which case $A$ is invertible in the usual sense) or $k=1$.

If $A$ is a polynomial in $P$ and $Q$, the functions $\det\Phi, \tr\Phi$ are also polynomial.
This allows us to simplify \eqref{dra} accordingly. To illustrate
things, consider a linear combination $A=aP+bQ$. In that case
\[ \Phi(t)=\begin{bmatrix} a+bt & b\sqrt{t(1-t)} \\ b\sqrt{t(1-t)} & b(1-t)\end{bmatrix}, \]
implying $\det\Phi(t)=ab(1-t)$ and $\tr\Phi(t)=a+b$. So, this particular $A$ is Drazin invertible
if and only if
$a=0$ or $b=0$ or $1\notin\sigma(H)$. Indeed,
if $a=b=0$, then $A=0$ is Drazin invertible.
If $a=0$ and $b$ is different from $0$, then $\Delta_2$ is empty
and the trace is separated from zero, so (\ref{dra}) holds.
Analogously for $b=0$ and $a$ different from zero.
Finally, if $ab$ is different from zero, then $\Delta_2$ is the whole spectrum with $1$ deleted,
and in order for $\det \Phi$ to be separated from zero
on it it is necessary and sufficient that the spectrum is separated from the point $1$.
But this is exactly the condition that $1$ is not in $\sigma(H)$.
Note that in all these cases, $A$ is also Moore-Penrose invertible.
We remark that the differences $P-Q$, along with some other simple polynomials in $P,Q$,
were treated by Deng \cite{Deng07}, prompting the considerations of \cite{BSpit09}.

\section{Compatible ranges} \label{s:cora}

As in \cite{Djik}, we will say that an operator $A$ acting on $\mathcal H$ has the
{\em compatible range} (CoR) {\em property} if $A$ and $A^*$ coincide
on $(\Ker A+\Ker A^*)^\perp$. It is easy to see (and was also observed in \cite{Djik})
that all the products $P,PQ,PQP,\ldots$ have this property. Those containing
an odd number of factors are Hermitian, which of course implies CoR. On the other hand,
the product of $n=2k$ interlacing $P$s and $Q$s is $A=(PQ)^k$.
So, $\Ker A\supset\Ker Q$, $\Ker A^*\supset\Ker P$, and in the notation of \eqref{canrep}
we have $(\Ker A+\Ker A^*)^\perp={\mathcal M}_{00}$. It remains to observe that
the restrictions of both $A$ and $A^*$ to this subspace are equal to the identity operator.

A somewhat tedious but straightforward computation of \[(\Ker A+\Ker A^*)^\perp\] with
the use of \eqref{ker} and its analogue for $A^*$ leads to the CoR criterion for arbitrary $A\in{\mathcal A}(P,Q)$
obtained in \cite{Spit181}. Namely:

The operator \eqref{elal} has the  CoR property if and only if $a_{ij}\in\R$ when
${\mathcal M}_{ij}\neq\{0\}$ and for (almost) every $t\in\Delta$ the matrix $\Phi(t)$ is either (i) Hermitian or (ii)
singular but not normal.

\section{A distance formula} \label{s:dist}

Along with $P,Q$, let us introduce the involution $U=2Q-I$. If $R$ is an orthogonal projection,
then, following \cite{Walt17}, $UR$ is called the symmetry of $R$
(with respect to $U$) and $R$ is said to be {\em orthogonal to its symmetry} if $RUR=0$.
Denote by ${\mathcal Q}_U$ the set of all orthogonal projections $R$
satisfying the orthogonality equation $RUR=0$.

It was shown in \cite{Walt17} that if $P$ is ``nearly orthogonal to its symmetry''
(quantitatively,
$x:=\norm{PUP}<\xi\approx 0.455$), then
\eq{dW} \dist(P,{\mathcal Q}_U)\leq \frac{1}{2}x+4x^2. \en
In fact, concentrating on ${\mathcal Q}^0_U:={\mathcal Q}_U\cap{\mathcal A}(P,Q)$ and
computing the norms along the lines of Section~\ref{s:rou} we
arrive at the following result established in~\cite{Spit18}:
if  ${\mathcal M}_{00}={\mathcal M}_{01}=\{0\}$ in \eqref{deco}, then
\eq{dist} \dist(P,{\mathcal Q}^0_U)= \sqrt{\frac{1}{2}\left(1-\sqrt{1-x^2}\right)}
=\frac{1}{2}x+\frac{1}{16}x^3+\cdots,  \en
and $ \dist(P,{\mathcal Q}^0_U)=1$ otherwise. Note that the latter case is
only possible if $\norm{PUP}=1$ and note also that there
are no a priori restrictions on $\norm{PUP}$ in order for \eqref{dist} to hold.

The distance \eqref{dist} is actually attained and, if $\norm{PUP}<1$,
the respective element of ${\mathcal Q}^0_U$ lies in ${\mathcal B}(P,Q)$.

\section{Index and trace} \label{s:index}

According to \cite{AvSeSiJFA}, $(P,Q)$ is a {\em Fredholm pair} if the operator
\eq{C} C:=QP\colon\im P\longrightarrow\im Q \en
is Fredholm, and the index $\ind(P,Q)$ of the pair $(P,Q)$ is by definition
the index of $C$. Using \eqref{M} and \eqref{imQ}, we can rewrite \eqref{C} in a more detailed form:
\[ C \colon {\mathcal M}_{00}\oplus {\mathcal M}_{01}\oplus {\mathcal M}\longrightarrow
{\mathcal M}_{00}\oplus {\mathcal M}_{10}\oplus{\mathcal N},\]
where ${\mathcal N}=\begin{bmatrix}\sqrt{H}\\ \sqrt{I-H}\end{bmatrix}({\mathcal M})$.

Now observe that $C$ acts as the identity on ${\mathcal M}_{00}$, the zero on ${\mathcal M}_{01}$,
while its action on ${\mathcal M}$ is the composition of the unitary operator
\[ \begin{bmatrix}\sqrt{H} \\ \sqrt{I-H}]\end{bmatrix}\colon {\mathcal M}\rightarrow {\mathcal N}  \]
with $\diag[\sqrt{H},\sqrt{H}]$.  We conclude that $\Ker C={\mathcal M}_{01}$
while $\im C$ is  the orthogonal sum of ${\mathcal M}_{00}$ with a dense
subspace of $\mathcal N$ which is closed if and only if the operator $H$ is invertible.
In particular, $(\im C)^\perp={\mathcal M}_{10}$.

So, the pair $(P,Q)$ is Fredholm if and only if ${\mathcal M}_{01},{\mathcal M}_{10}$
are finite-dimensional and $H$ is invertible. Moreover, if these conditions hold, then
\[ \ind(P,Q)= \dim{\mathcal M}_{01}-\dim{\mathcal M}_{10}.\]
This result can be recast in terms of the difference $P-Q$. Namely, the operator $H$ is
invertible if and only if $\pm 1$ are at most isolated points of $\sigma(P-Q)$ (see formula \eqref{spq} and
the explanations following it), while ${\mathcal M}_{01}$ and ${\mathcal M}_{10}$
are simply the eigenspaces of $P-Q$ corresponding to $\pm 1$, due to \eqref{p-q}. We thus
arrive at Proposition 3.1 of \cite{AvSeSiJFA}, which says that
the pair $(P,Q)$ is Fredholm if and only if $\pm 1$ are (at most) isolated points of
$\sigma(P-Q)$ having finite multiplicity and that under these conditions
\eq{ind} \ind (P,Q)=\dim\Ker(P-Q-I)-\dim\Ker(P-Q+I). \en
Because $\sigma(P-Q)\subset[-1,1]$, we see in particular that if $P,Q$ are in generic position, then the pair is
Fredholm if and only if \eq{sno} \norm{P-Q}<1,\en
and then $\ind(P,Q)=0$. This was pointed out in \cite{AmSi}.

Let us now consider powers of $P-Q$. Since \eqref{fpq} may be rewritten as
\[ \Phi_{P-Q}=\sqrt{1-t}
\begin{bmatrix} \sqrt{1-t} & -\sqrt{t} \\  -\sqrt{t} & - \sqrt{1-t}\end{bmatrix}, \]
with the matrix factor on the right-hand side being an involution, it is easy to see
that \eqref{p-q} implies that, for every even $k=2n$,
\[ (P-Q)^k=0_{{\mathcal M}_{00}}\oplus I_{{\mathcal M}_{01}}\oplus I_{{\mathcal M}_{10}}\oplus
0_{{\mathcal M}_{11}}\oplus \diag[(I-H)^n,(I-H)^n]. \]
Consequently, for odd powers $k=2n+1$,
\begin{multline}\label{ppq} (P-Q)^k=0_{{\mathcal M}_{00}}\oplus I_{{\mathcal M}_{01}}\oplus
(-I)_{{\mathcal M}_{10}}\oplus 0_{{\mathcal M}_{11}}
\\ \oplus (I-H)^{n+1/2}\begin{bmatrix} \sqrt{I-H} & -\sqrt{H} \\ -\sqrt{H} & -\sqrt{I-H}\end{bmatrix}. \end{multline}
Suppose now that for some $m$  the $m$-th power of $P-Q$ is a trace class operator.
Then ${\mathcal M}_{01}, {\mathcal M}_{10}$ are finite-dimensional, and for every $k\geq m$
the last direct summand in
\eqref{ppq} is a zero-trace operator. We thus have
\eq{tr} \tr(P-Q)^k = \dim\Ker(P-Q-I)-\dim\Ker(P-Q+I) \en
independently of $k$.

Note also that $(P-Q)^k$ being a trace class operator implies that $P-Q$, and therefore $I-H$, is compact.
Then $H$, as a Fredholm operator with zero index and (by its construction) satisfying
$\Ker H=\{0\}$ is in fact invertible. As stated above, the pair $(P,Q)$ is thus Fredholm,
and \eqref{ind} holds. Comparing \eqref{ind} with \eqref{tr}, we arrive at the formula
\[ \tr(P-Q)^k =  \ind (P,Q) \] valid for any odd $k\geq m$ provided that $(P-Q)^m$ is trace class.
This is \cite[Theorem 4.1]{AvSeSiJFA}. In relation to their
physics applications, the results of this section are also treated in \cite{AvSeSiJMP}.

\section{Intertwining} \label{s:inter}

In the early 1950s, Kato (unpublished) found a unitary operator $U$ satisfying $UP=QU$
provided that \eqref{sno} holds. In \cite{AvSeSiJFA} it was established that, under the same
condition \eqref{sno}, the unitary $U$ can be constructed to satisfy the two equations
\eq{UPQ} UP=QU \text{ and } UQ=PU; \en
we will say that such a $U$ {\em intertwines} $P$ with $Q$.

A necessary and sufficient condition for such $U$ to exist
is that in \eqref{deco} \eq{010}  \dim{\mathcal M}_{01}=\dim{\mathcal M}_{10}, \en
see \cite[Theorem 6]{WangDuDou}.  Note that \eqref{sno} implies  \eq{000} {\mathcal M}_{01}={\mathcal M}_{10}=\{0\}, \en
so that \eqref{010} holds in a trivial way.

A description of all $U$ satisfying \eqref{UPQ} was provided in \cite{DSCD}.
In the notation \eqref{canrep} it looks as follows \cite{BoSiS}:
\eq{U} U=U_0 \oplus \begin{bmatrix} 0 & U_{10} \\ U_{01} & 0\end{bmatrix}\oplus U_1
\oplus{\mathbf W}^*\begin{bmatrix} V & 0 \\ 0 & V \end{bmatrix}
\begin{bmatrix}  \sqrt{H} & \sqrt{I-H} \:\\ \sqrt{I-H} & -\sqrt{H} \end{bmatrix}{\mathbf W}. \en
Here $U_j$, $U_{ij}$ are arbitrary unitary operators acting on ${\mathcal M}_{jj}$
and from ${\mathcal M}_{ji}$ onto ${\mathcal M}_{ij}$, respectively, and
$V$ is an arbitrary unitary operator acting on $\mathcal M$ and commuting with $H$.

Invoking \eqref{elal}, it was also observed in \cite{BoSiS}  that operators $U$ intertwining
$P$ and $Q$ can be chosen in  ${\mathcal A}(P,Q)$ only if
instead of \eqref{010} the stronger condition \eqref{000} is imposed.
All such operators $U$ are then given by
\eq{UA}
U=a_0I_{{\mathcal M}_{00}}\oplus a_1I_{{\mathcal M}_{11}}\oplus
{\mathbf W}^*\begin{bmatrix} \phi(H) & 0 \\ 0 & \phi(H) \end{bmatrix}
\begin{bmatrix}  \sqrt{H} & \sqrt{I-H}\: \\ \sqrt{I-H} & -\sqrt{H} \end{bmatrix}{\mathbf W},
\en
where $\abs{a_0}=\abs{a_1}=1$ and $\phi$ is a Borel measurable unimodular function on $[0,1]$.

In its turn, such $U$ lie in ${\mathcal B}(P.Q)$ if and only if the unimodular function
$\phi$ is continuous on $[0,1]$, not just measurable.
Finally, if the pair $P,Q$ is in  generic position  and the spectrum of $PQP$ is simple,
then all operators satisfying \eqref{UPQ} lie in ${\mathcal A}(P,Q)$.

\section{The supersymmetric approach} \label{s: super}

The pertinent results of Sections~\ref{s:index} and \ref{s:inter} were
obtained in \cite{AvSeSiJMP,AvSeSiJFA} solely based on
the simple (and directly verifiable)
observation that for any two orthogonal  projections $P,Q$ the (selfadjoint) operators
\eq{AB} A=P-Q,\quad B=I-P-Q \en satisfy
\eq{SS} A^2+B^2=I, \quad AB+BA=0. \en
Because of the second formula in \eqref{SS}, it is natural to speak of the supersymmetric approach.

The approach of \cite{AmSi,DouDuWang17,DSCD} and \cite{WangDuDou} was geometrical,
using either \eqref{canrep} or its equivalents.
In \cite{Simon17} a point was made to derive the existence criterion for the
intertwining unitary $U$ via the supersymmetric approach.
For the description of all such $U$, this was done in \cite[Section 4]{BoSiS}.
Here we would like to show how the reasoning of the latter, with
some modifications, can be used to derive Halmos' canonical representation \eqref{deco},\eqref{canrep}  for
\eq{PQAB} P=\frac{1}{2}(I+A-B), \quad  Q=\frac{1}{2}(I-A-B) \en
directly from \eqref{SS}.

The first of the formulas \eqref{SS} shows that  the restriction of $B$ to the
eigenspaces of $A$ corresponding to the eigenvalues $\pm 1$ equals zero.
Denote these eigenspaces by ${\mathcal M}_{01}$ and ${\mathcal M}_{10}$, respectively,
and let $\langle \cdot, \cdot \rangle$ be the scalar product in
${\mathcal H}$. Then for a unit vector $x\in{\mathcal M}_{01}$
we have $\scal{(P-Q)x,x}=1$. But both $\scal{Px,x}$ and $\scal{Qx,x}$ take their values
in $[0,1]$, leaving us with the only option $\scal{Px,x}=1$, $\scal{Qx,x}=0$.
This in turn implies $Px=x$ and $Qx=0$, i.e., $P|{\mathcal M}_{01}=I$, $Q|{\mathcal M}_{01}=0$.
Similarly, $P|{\mathcal M}_{10}=0$, $Q|{\mathcal M}_{10}=I$.
This agrees with \eqref{canrep} and allows us to consider now the restrictions of $A,B$ to the
orthogonal complement ${\mathcal H}'$ of ${\mathcal M}_{01}\oplus{\mathcal M}_{10}$.

Representing ${\mathcal H}'$ as the orthogonal sum of $\Ker A$ and the spectral subspaces
of $A$ corresponding to the positive (resp., negative) parts of its spectrum,
we can write $A|{\mathcal H}'$ and $B|{\mathcal H}'$ as $A'=\diag[0,A_+,-A_-]$ and
\[ B'=\begin{bmatrix} B_{00} & B_{01} & B_{02} \\ B_{01}^* & B_{11} & B_{12} \\ B_{02}^* & B_{12}^* & B_{22}\end{bmatrix}, \]
with $A_\pm$ being positive definite operators.

Since formulas \eqref{SS} carry over to $A',B'$, we have in particular \[ A_+B_{11}+B_{11}A_+=0.\]
Thus, the operator $A_+B_{11}$ has zero Hermitian part,
and so its spectrum is purely imaginary. Since for any two operators $X,Y$
we have $\sigma(XY)\cup\{0\}=\sigma(YX)\cup\{0\}$, the spectrum of
$A_+^{1/2}B_{11}A^{1/2}$ also is purely imaginary. On the other hand, the latter
operator is selfadjoint, and hence its spectrum is real. Combining these
two observations we conclude that the selfadjoint operator $A_+^{1/2}B_{11}A^{1/2}$
has zero spectrum and thus itself is zero. From the injectivity of $A_+$
we conclude that $B_{11}=0$.
Similarly, \[A_-B_{22}+B_{22}A_-=0\] implies that $B_{22}=0$.

With these simplifications in mind, the second part of \eqref{SS} is now equivalent to
\eq{BA1} B_{01}A_+=0, \quad B_{02}A_-=0,\en and \eq{BA2} A_+B_{12}=B_{12}A_-. \en
Invoking the injectivity of $A_\pm$ again, we see from \eqref{BA1} that the blocks
$B_{01},B_{02}$ are also equal to zero, and so $B'$ takes the form
\[ B_{00}\oplus  \begin{bmatrix} 0 & B_{12} \\  B_{12}^* & 0\end{bmatrix}. \]
In particular, $\Ker A$ is an invariant subspace of $B$. According to the first formula in \eqref{SS},
the restriction $B_{00}$ of $B$ to $\Ker A$ is a (selfadjoint) involution.
Consequently, $\Ker A$ splits into the orthogonal sum of the eigenspaces of $B$ corresponding
to the eigenvalues $\pm 1$. Denoting them by
${\mathcal M}_{00}$ and ${\mathcal M}_{11}$ and using \eqref{PQAB},
we find ourselves in agreement with \eqref{canrep} again.

With a slight abuse of notation, we are now left with the following. Let $A,B$ be given by
\[ A=\begin{bmatrix} A_+ & 0 \\ 0 & -A_-\end{bmatrix}, \quad B=
\begin{bmatrix} 0 & B_{12} \\  B_{12}^* & 0\end{bmatrix}, \]
with $A_\pm$ positive definite and not having $1$ as an eigenvalue.
Let also \eqref{BA2} hold and suppose
\eq{AB2} A_+^2+B_{12}B_{12}^*=I, \quad A_-^2+B_{12}^*B_{12}=I. \en
Our task is to show that then the pair \eqref{PQAB} is in generic position
and admits the respective representation \eqref{canrep}.

Of course, \eqref{AB2} is simply the first condition in \eqref{SS} written block-wise.

Since $1$ is not an eigenvalue of $A_\pm$, equalities \eqref{AB2} imply that $B_{12}$ has
zero kernel and dense range.
In its polar representation \[ B_{12}=CV, \quad C=\sqrt{B_{12}^*B_{12}} \]
the operator $V$ is an isometry between the domains of $A_\pm$, implying in particular
that these domains have equal dimensions.
The unitary similarity $\diag[I,V]$ allows us to replace the pair $(A,B)$ by
\[ \begin{bmatrix} A_+ & 0 \\ 0 & -VA_-V^*\end{bmatrix}, \quad B= \begin{bmatrix} 0 & C \\ C & 0\end{bmatrix},\]
for which \eqref{AB2} turns into $A_+^2+C^2= (VA_-V^*)^2+C^2=I$.

But $A_+$ and $VA_-V^*$ are both positive definite. So, the latter equality defines them uniquely as
\[ A_+=VA_-V^*=\sqrt{I-C^2}:=S.\] We have thus found a unitary similarity under which $P,Q$
become
\eq{PQm} P=\frac{1}{2}\begin{bmatrix} I+S & -C \\ -C &   I-S\end{bmatrix}, \quad
Q=\frac{1}{2}\begin{bmatrix} I-S & -C \\ -C &   I+S\end{bmatrix}. \en
A side note: the representation \eqref{PQm}, being more ``balanced'', has some advantages
over the generic portion of \eqref{canrep}.
In particular, the existence of an intertwining $U$ becomes obvious:
the permutation $\begin{bmatrix} 0 & I \\ I & 0\end{bmatrix}$ does the job.

For the task at hand, however, one more unitary similarity is needed, one which
reduces $P$ from \eqref{PQm} to the form $\diag[I,0]$.
To this end, let us introduce the selfadjoint involution
\[ J=\frac{\sqrt{2}}{2}(I+S)^{-1/2}\begin{bmatrix}C & -(I+S) \\  -(I+S) & -C\end{bmatrix}. \]
A direct computation shows that then indeed
$JPJ=\diag[I,0]$, while
\[ JQJ = \begin{bmatrix} S^2 & CS \\ CS & C^2\end{bmatrix}. \]
It remains to relabel $C^2=H$.

So, it is not surprising that any result pertinent to pairs of orthogonal projections
can be derived from scratch just by using
the purely algebraic relations \eqref{SS}. The supersymmetric approach is an Answerer that can rival
with Halmos' theorem.

\bigskip
\noindent
Albrecht B\"ottcher\\
Fakult\"at f\"ur Mathematik\\
TU Chemnitz\\ D-09107 Chemnitz\\ Germany\\
{\tt aboettch@math.tu-chemnitz.de}

\bigskip
\noindent
Ilya M. Spitkovsky\\
Division of Science\\  New York  University Abu Dhabi (NYUAD)\\ Saadiyat Island\\
P.O. Box 129188 Abu Dhabi\\ UAE\\
{\tt ims2@nyu.edu, imspitkovsky@gmail.com}

\end{document}